\documentclass[a4paper]{article}
\usepackage[utf8]{inputenc}
\usepackage[ngerman, english]{babel}
\usepackage{amsmath,amsthm,amssymb,amsfonts}
\usepackage{mathrsfs}
\usepackage{ifthen}
\usepackage{enumerate}
\usepackage{comment}
\usepackage{subcaption}
\usepackage{hyphsubst}
\usepackage{xcolor}
\usepackage{calc}
\usepackage{graphicx}
\usepackage[toc,page]{appendix}

\makeatletter
\def\@hspace#1{\begingroup\setlength\dimen@{#1}\hskip\dimen@\endgroup}
\makeatother

\newtheorem{theorem}{Theorem}[section]
\newtheorem{lemma}[theorem]{Lemma}
\newtheorem{definition}[theorem]{Definition}

\newtheorem{assumption}[theorem]{Assumption}

\newcommand{\dd}{\mathrm{d}}
\newcommand{\nv}{\nu}

\newcommand{\ol}[1]{\overline{#1}}

\newcommand{\spl}{\langle}
\newcommand{\spr}{\rangle}
\newcommand{\bpm}{\begin{pmatrix}}
\newcommand{\epm}{\end{pmatrix}}

\DeclareMathOperator{\diag}{diag}

\renewcommand{\dim}{\operatorname{dim}}

\renewcommand{\div}{\operatorname{div}}

\DeclareMathOperator{\supp}{supp}

\newcommand{\setC}{\mathbb{C}}

\newcommand{\setN}{\mathbb{N}}

\newcommand{\setR}{\mathbb{R}}

\definecolor{brickred}{rgb}{0.8, 0.25, 0.33}
\definecolor{bostonuniversityred}{rgb}{0.8, 0.0, 0.0}
\definecolor{cornellred}{rgb}{0.7, 0.11, 0.11}
\definecolor{corn}{rgb}{0.98, 0.93, 0.36}
\definecolor{schoolbusyellow}{rgb}{1.0, 0.85, 0.0}
\definecolor{TUblue}{rgb}{0,102,153}
\colorlet{TUbluelight}{TUblue!30!white}

\title{Radial complex scaling for anisotropic scalar resonance problems}
\author{Martin Halla\footnote{Max-Planck-Institut f\"ur Sonnensystemforschung, Justus-von-Liebig-Weg 3, 37077 G\"ottingen, Deutschland (halla@mps.mpg.de)} \footnote{Institut f\"ur Numerische und Angewandte Mathematik, Georg-August Universit\"at G\"ottingen, Lotzestraße 16-18, 37083 Göttingen, Deutschland}}

\begin{document}

\maketitle
\begin{abstract}
\noindent
The complex scaling/perfectly matched layer method is a widely spread technique to simulate wave propagation problems in open domains.
The method is very popular, because its implementation is very easy and does not require the knowledge of a fundamental solution.
However, for anisotropic media the method may yield an unphysical radiation condition and lead to erroneous and unstable results.
In this article we argue that a radial scaling (opposed to a cartesian scaling) does not suffer from this drawback and produces the desired radiation condition.
This result is of great importance as it rehabilitates the application of the complex scaling method for anisotropic media.
To present further details we consider the radial complex scaling method for scalar anisotropic resonance problems.
We prove that the associated operator is Fredholm and show the convergence of approximations generated by simulateneous domain truncation and finite element discretization.
We present computational studies to undergird our theoretical results.
\end{abstract}

\section{Introduction}\label{sec:introduction}

Since the 1970s the complex scaling (CS) method has been used in molecular physics to study resonances \cite{Simon:78,Moiseyev:98}.
In the 1990s B\'erenger introduced his perfectly matched layer (PML) method \cite{Berenger:94} as reflectionless sponge layer for time-dependent electromagnetic wave equations, and
this method became popular for all kinds of wave propagation problems.
Soon the PML method was recognized to be a complex scaling technique \cite{ChewWeedon:94,CollinoMonk:98a}.
Also the original variant of the Hardy space infinite element method \cite{HohageNannen:09,Halla:16} was recognized in \cite{NannenWess:19} to be a conjunction of a complex scaling and an infinite element method.
The great advantage of the CS/PML method is that its application does not require the knowledge of a fundamental solution, and hence the method is very easy to implement.
In particular, for time-harmonic equations the method requires only the implementation of some complex valued coefficient functions, and
for time-dependent equations the method requires only the introduction of some auxiliary fields.
No doubt, this versality lead to the tremendous popularity of the CS/PML method.

Even though the method is easy to implement, its naive application does not always yield meaningful results.
For certain types of equations convenient complex scalings do not yield physically correct and stable CS/PML methods, and it is necessary to construct problem adapted scalings.
To be physically meaningful the complex scaling has to be designed such that evanescent waves stay evanescent and propagative waves with positive group velocity become evanescent.
In particular, if an equation admits backward waves, i.e.\ waves for which the product of the phase velocity vector and the group velocity vector is negative, convenient CS/PML methods do not work.
Backward waves can occur e.g. for dispersive equations \cite{CassierJolyKachanovska:17,BecacheKachanovska:17}.
Typically for such equations either all waves are backward or all waves are forward, depending on the frequency.
Differently, in waveguides it is possible that for some frequencies forward and backward waves exist simultaneously \cite{BonnetBDChambeyronLegendre:14,HallaNannen:15,HallaHohageNannenSchoeberl:16}.

An other issue is raised for anisotropic \cite{BecacheFauqueuxJoly:03,BecFliKacKaz:21} and advective equations \cite{BecacheBonnetBDLegendre:06}.
A main tool in the analysis of such equations in homogeneous open domains has been the characterization of materials in terms of their slowness curves \cite{BecacheFauqueuxJoly:03}.
If the slowness curves are concentric ellipsoids, then the axes of a cartesian complex scaling can be aligned to the axes of the slowness curves and such methods seem to produce reliable results.
Nevertheless, for certain media (e.g.\ some orthotropic elastic materials) this is not possible and the CS/PML method is reported to fail \cite{BecacheFauqueuxJoly:03}.
Thus the construction of CS/PML methods for ``the case of anisotropic elastic waves, remains (\dots) a challenging open problem'' \cite{Joly:12}.
Note that this phenomenon is seperate from the occurence of backward waves, and for ``convenient'' equations as considered in \cite{BecacheFauqueuxJoly:03,Joly:12} there exist no backward waves.

These drawbacks of CS/PML methods gave rise to the development of novel alternative methods such as the two-pole Hardy space infinite element method \cite{HallaNannen:15,HallaHohageNannenSchoeberl:16,HallaNannen:18} and the halfspace matching method \cite{BonnetBDFlissTonnoir:18,BonnetBDetal:20}.

However, in this article we argue that \emph{radial} complex scalings can be applied succesfully to anisotropic equations.
The idea is actually very simple.
If the asymptotic behavior of outgoing solutions is governed by $e^{i\omega rf(\hat x)}/r$ with a uniformly positive function $f$ and polar coordinates $x=r\hat x$, then a radial scaling $\tilde r(r)=(1+\frac{\alpha_0}{-i\omega})r$ with $\alpha_0>0$ achieves the damping of outgoing waves:
\begin{align*}
\left|\frac{e^{i\omega \tilde r(r)f(\hat x)}}{\tilde r(r)}\right|\leq \frac{e^{-\alpha_0 \inf_{\hat x\in S^2}\{f(\hat x)\} r}}{\sqrt{1+\alpha_0^2/\omega^2}r}.
\end{align*}
This very basic observation has the powerful consequence that it rehabilitates the application of CS/PML methods for anisotropic equations.
For scalar and certain elastic equations \cite[Appendix]{Vavrycuk:07} explicit expressions of Greens functions are known and these Green functions admit the mentioned behavior.

An additional argument can be made for time-dependent equations.
Actual stability and convergence results on CS/PML methods in the time-domain are rare \cite{DiazJoly:06,BecacheKachanovska:20}.
Most reports deal only with constant coefficient scalings and such an analysis can already be challenging.
However, for radial constant coefficient scalings the obtained equation is very similar to the one-dimensional case, which can be analyzed quite easily \cite{Joly:12}.
Consider e.g. an anisotropic elastic inital value problem
\begin{align*}
\partial_{tt}u-\div \sigma(u)=0 \quad\text{ in }\setR^3,  \quad u(x,t=0)=u_0(x),
\end{align*}
with displacement vector $u$, strain tensor $\epsilon(u):=1/2(\nabla u+\nabla u^\top)$, stress tensor $\sigma(u)=C\epsilon(u)$ and fourth-oder tensor $C$.
Then the radial constant coefficient scaling $\tilde r(r)=(1+\frac{\alpha_0}{-i\omega})r$ yields the equation
\begin{align*}
\partial_{tt}\tilde u+2\alpha_0\partial_t\tilde u+\alpha_0^2\tilde u-\div \sigma(u)=0 \quad\text{ in }\setR^3,  \quad \tilde u(x,t=0)=\tilde u_0(x).
\end{align*}
Thus for a solution $\tilde u$ we can test with $\partial_t\tilde u$ and derive an equation for the energy
\begin{align*}
E(t)&:=\|\partial_t \tilde u(\cdot,t)\|^2_{L^2(\setR^3)}+\|\alpha_0\tilde u(\cdot,t)\|^2_{L^2(\setR^3)}+\|C^{1/2}\epsilon(\tilde u(\cdot,t))\|^2_{L^2(\setR^3)}, \\
\partial_t E(t)&=-2\alpha_0 \|\partial_t \tilde u(\cdot,t)\|^2_{L^2(\setR^3)}\leq 0.
\end{align*}
Hence the energy $E(\tilde u)$ is bounded in $t\geq0$ and thus the equation is stable.
For this reasoning it is completely irrelevant, if the material described by the tensor $C$ is isotropic or anisotropic.

Although, the correct damping of outgoing waves and stability for constant coefficient scalings is only the first step in the analysis and does not predict any Fredholm or stability results.
To investigate these questions we start in this article with the analysis of the ``most simple'' case: the resonance problem for the anisotropic scalar wave equation.
That is we consider a symmetric positive definite matrix $\varsigma\in\setR^{3\times 3}$ and search for solutions $(\omega,u)$ to
\begin{align*}
-\div(\varsigma\nabla u)-\omega^2u =0 \quad+\text{boundary and radiation conditions}.
\end{align*}
We plan to extend our results to anisotropic elastic and time-dependent equations in forthcoming articles.
For the anisotropic scalar wave equation it turns out that the framework for the isotropic case \cite{Halla:19Diss,Halla:20PML} can be reused and the analysis requires only slight modifications.
In particular, in Lemma~\ref{lem:TildeUofU} we obtain that the radial complex scaling produces the desired damping of outgoing waves.
In Theorem~\ref{thm:wTc} we report the weak T-coercivity (equivalent to Fredholmness with index zero) of the operator corresponding to the eigenvalue problem.
The proof requires an additional Lemma~\ref{lem:NumRan} on the numerical range of certain complex valued matrices.
Lemma~\ref{lem:NumRan} gives rise to an extra assumption which involves the complex scaling parameters and the degree of anisotropicy of $\varsigma$.
This assumption is important for time-depent problems, because it shows that the convenient frequency dependency $1/(-i\omega)$ of the scaling is likely to produce instabilities.
As a remedy we propose a frequency dependency $1/(c-i\omega)$ with a positive constant $c$.
See the discussion at the end of Section~\ref{sec:wTc}.
The proof of convergence of approximations by means of simultaneous domain truncation and discretization is actual identical to the isotropic case.

The remainder of this article is structured as follows.
In Section~\ref{sec:resonanceproblem} we introduce the resonance problem and its CS/PML reformulation.
In Lemmata~\ref{lem:TildeUofU} and~\ref{lem:TildeUsolution} we report the damping of outgoing waves.
In Section~\ref{sec:wTc} we proof the weak T-coercivity of the associated operator.
In Section~\ref{sec:Appr} we proof the convergence of approximations.
In Section~\ref{sec:comp} we present computational studies to undergird our theoretical results.
We close the article with a conclusion in Section~\ref{sec:conclusion}.

\section{The resonance problem and its reformulation}\label{sec:resonanceproblem}

Let $B_r(x_0)\subset\setR^3$ be the open ball with radius $r>0$ centered at $x_0\in\setR^3$ and $B_r:=B_r(0)$.
For a Lipschitz domain $D\subset\setR^3$ let $\nv$ be its outward unit normal vector at $\partial D$ and
\begin{align*}
\tilde H^1_\mathrm{loc}(D)&:=\{u\colon u|_{D\cap B_r}\in H^1(D\cap B_r) \text{ for all }r>0\text{ with }D\cap B_r\neq\emptyset\},\\
\tilde H^1_\mathrm{0,loc}(D)&:=\{u\in \tilde H^1_\mathrm{loc}(D)\colon u|_{\partial D}=0\}.
\end{align*}
Let $\Omega\subset\setR^3$ be a connected Lipschitz domain such that the complement $\Omega^c$ is compact and non-empty.
Let $\varsigma\in\setR^{3\times3}$ be symmetric positive definite, and $\varsigma_\mathrm{min}$ and $\varsigma_\mathrm{max}$ be its minimal and maximal eigenvalues.
Then we seek solutions $(\omega,u)$ to
\begin{subequations}\label{eq:ResonanceProblemFormal}
\begin{align}
\label{eq:ResonanceProblemFormalA}
-\div \varsigma \nabla u-\omega^2 u&=0\quad\text{in }\Omega,\\
\label{eq:ResonanceProblemFormalB}
u&=0\quad\text{on }\partial\Omega,
\end{align}
together with the abstract radiation condition (which will be specified in Definition~\ref{def:RadiationCondition})
\begin{align}
\label{eq:ResonanceProblemFormalC}
u&\text{ is outgoing}.
\end{align}
\end{subequations}
Let $r_0>0$ be such that $\Omega^c\subset B_{r_0}$.
Let $\supp u:=\ol{\{x\in\Omega\colon u(x)\neq0\}}^\mathrm{cl}$ be the support of a function $u$.
Note that for inhomogeneous materials $\hat\varsigma(x)$ with density $\rho(x)$,  and/or scattering problems with right hand side $f\in L^2(\Omega)$ the radius $r_0$ needs to be big enough such that $\supp f, \supp (\rho-1), \supp (\hat\varsigma-\varsigma) \subset B_{r_0}$.
The physical meaningful radiation condition for the isotropic wave equation is that $u$ admits a representation
\begin{align*}
u(x)=\frac{i\omega}{4\pi} \int_{\partial B_{r_0}}
u(y) \nabla_y G_\omega(x,y)\cdot\nv(y) - G_\omega(x,y) \nabla_y u(y)\cdot\nv(y) \dd y, \quad x\in B_{r_0}^c
\end{align*}
with the fundamental solution $G_\omega(x,y):=h^{(1)}_0(\omega|x-y|)=\frac{e^{i\omega |x-y|}}{i\omega |x-y|}$.
By means of the transformation rule we derive an analogous radiation condition for the anisotropic case.
\begin{definition}[Radiation condition]\label{def:RadiationCondition}
Let $(\omega,u)$ be a solution to~\eqref{eq:ResonanceProblemFormal}.
We call $u$ to be outgoing if it admits a representation
\begin{align*}
u(x)=\frac{i\omega}{4\pi} \int_{\partial B_{r_0}}
u(y) \nabla_y G_{\omega,\varsigma}(x,y)\cdot\nv(y) - G_{\omega,\varsigma}(x,y) \nabla_y u(y)\cdot\nv(y) \dd y, \,\,\,\,\, x\in B_{r_0}^c
\end{align*}
with the fundamental solution
$G_{\omega,\varsigma}(x,y):=\det\varsigma^{-1/2} h^{(1)}_0(\omega|\varsigma^{-1/2}(x-y)|)$.
\end{definition}

\subsection{Radial complex scaling}

We will define a complex change of the radial coordinate
$\tilde r(r)=(1+\gamma\tilde\alpha(r))r$ in terms of a profile function $\tilde\alpha$ and a constant $\gamma$.
We make assumptions on those as follows.
\begin{assumption}\label{ass:TildeAlpha}
Let $r_1>r_0$, $\gamma\in\{z\in\setC\colon\Re(z)\geq0, \Im(z)>0\}$,
and $\tilde\alpha\colon\setR^+_0\to\setR^+_0$ be such that
\begin{enumerate}
 \item\label{it:talpha0} $\tilde\alpha(r)=0$ for $r\leq r_1$,
 \item\label{it:talphacont} $\tilde\alpha$ is continuous,
 \item\label{it:nontrivial} $\tilde\alpha(r)>0$ for $r>r_1$,
 \item\label{it:talphamonoton} $\tilde\alpha$ is non-decreasing,
 \item\label{it:talphaC2} $\tilde\alpha$ is twice continuously differentiable in $(r_1,+\infty)$ with continuous
 extensions of $\partial_r\tilde\alpha$ and $\partial_r\partial_r\tilde\alpha$ to $[r_1,+\infty)$,
 \item\label{it:taabounded} $\tilde\alpha$ and $\alpha$ are bounded.
\end{enumerate}
\end{assumption}
Note that the last item of Assumption~\ref{ass:TildeAlpha} is a luxury.
It allows us to avoid weigthed Sobolev spaces and thus simplifies the presentation.
See \cite{Halla:19Diss,Halla:20PML} for the treatment of unbounded profile functions.

In the following we introduce additional functions which will all depend on $\tilde\alpha$ and $\gamma$.
These auxiliary functions will be necessary to formulate the forthcoming theory.
Let
\begin{subequations}\label{eq:ComplexScalingQuantities}
\begin{align}
\label{eq:tdDef} \tilde d(r)&:=1+\gamma\tilde\alpha(r),\\
\label{eq:trDef} \tilde r(r)&:=\tilde d(r)r,\\
\label{eq:alphaDef} \alpha(r)&:=r\partial_r\tilde\alpha(r)+\tilde\alpha(r),\\
\label{eq:dDef} d(r)&:=1+\gamma\alpha(r),\\
\label{eq:dzDef} d_0&:=\displaystyle\lim_{r\to+\infty}(\tilde d(r)/|\tilde d(r)|).
\end{align}
The definitions of $\alpha$ and $d$ have to be understood piece-wise.
We note that the limit in~\eqref{eq:dzDef} exists in $\setC$ due to Assumption~\ref{ass:TildeAlpha}.
The function $d$ is chosen such that $\partial_r\tilde r(r)=d(r)$.
For $f=\tilde\alpha,\alpha,\tilde d,d,\tilde r$ we adopt the overloaded notation
\begin{align}\label{eq:OverloadedNotation}
f(x):=f(|x|), \quad x\in\setR^3.
\end{align}
\end{subequations}
Let the spherical coordinates $r>0, \phi\in[0,2\pi), \theta\in[0,\pi]$ be defined by
\begin{align*}
x&=r \hat x(\phi,\theta),
\qquad \hat x(\phi,\theta):=(\sin\theta\cos\phi,\sin\theta\sin\phi,\cos\theta)^\top.
\end{align*}
For a solution $u$ to \eqref{eq:ResonanceProblemFormal} consider the complex scaled function
\begin{align*}
\tilde u(x):=u(\tilde r(r)\hat x(\phi,\theta)).
\end{align*}
In order to guarantee that this definition is meaningful we introduce the next assumption and lemmata.

\begin{assumption}\label{ass:r1r0}
Let $r_0$, $\varsigma_\mathrm{min}$ and $\varsigma_\mathrm{max}$ be as defined previously in this section and $r_1$ be as in Assumption~\ref{ass:TildeAlpha}.
Let $r_1>\frac{\varsigma_\mathrm{max}}{\varsigma_\mathrm{min}}r_0$.
\end{assumption}
We require Assumption~\ref{ass:r1r0} for Lemma~\ref{lem:dxy} and subsequently for Lemma~\ref{lem:TildeUofU}.
However, it is possible that Assumption~\ref{ass:r1r0} is not necessary and is just an artefact of our technique applied.

\begin{lemma}\label{lem:dxy}
Let Assumptions~\ref{ass:TildeAlpha} and \ref{ass:r1r0} be satisfied.
Let $x\in\Omega$ and $y\in \partial B_{r_0}$.
Then $(\tilde r_x \hat x-r_y \hat y)^\top \varsigma^{-1} (\tilde r_x \hat x-r_y \hat y)$ $\in$
$\{z\in\setC\colon \Im(z)>0\}\cup \setR_{\geq0}$.
Further, it holds $(\tilde r_x \hat x-r_y \hat y)^\top \varsigma^{-1} (\tilde r_x \hat x-r_y \hat y)=0$ if and only if $x=y$.
\end{lemma}
\begin{proof}
We proceed as in~\cite[Lemma~4.1]{BermudezHervellaNPrietoRodriguez:08}.
For $|x|\leq r_1$ it holds $\tilde r_x(x)=|x|$ and hence the claim follows.
For $|x|>r_1$ we apply $\Re(\tilde r_x)\geq r_x$, $\Im(\tilde r_x)>0$, and compute by means of Young's inequality that
\begin{align*}
\Im\big( (\tilde r_x \hat x&-r_y \hat y)^\top \varsigma^{-1} (\tilde r_x \hat x-r_y \hat y) \big)\\
&=2\Re(\tilde r_x)\Im(\tilde r_x) \hat x^\top\varsigma^{-1}\hat x
-2r_y\Im(\tilde r_x) \hat x^\top \varsigma^{-1}\hat y\\
&\geq 2\Re(\tilde r_x)\Im(\tilde r_x) \hat x^\top\varsigma^{-1}\hat x
-r_y\Im(\tilde r_x) (\hat x^\top \varsigma^{-1}\hat x+\hat y^\top \varsigma^{-1}\hat y)\\
&=\Im(\tilde r_x) \Big( \hat x^\top \varsigma^{-1}\hat x(\Re(\tilde r_x)-r_y)
+\hat x^\top \varsigma^{-1}\hat x \Re(\tilde r_x) -\hat y^\top \varsigma^{-1}\hat y r_y\Big)\\
&\geq \Im(\tilde r_x) \Big( \varsigma^{-1}_\mathrm{max}(\Re(\tilde r_x)-r_y)
+\varsigma^{-1}_\mathrm{max} \Re(\tilde r_x) -\varsigma^{-1}_\mathrm{min} r_y\Big)\\
&\geq \Im(\tilde r_x) \Big( \varsigma^{-1}_\mathrm{max}(r_1-r_0)
+\varsigma^{-1}_\mathrm{max} r_1 -\varsigma^{-1}_\mathrm{min} r_0\Big).
\end{align*}
Hence, the claim follows.
\end{proof}
Due to Lemma~\ref{lem:dxy} we can define the function
\begin{align*}
d_\varsigma(x,y):=\sqrt{(\tilde r_x \hat x-r_y \hat y)^\top \varsigma^{-1} (\tilde r_x \hat x-r_y \hat y)}, \quad x\in\Omega, y\in\partial B_{r_0}
\end{align*}
such that $\Im(d_\varsigma(x,y))\geq0$.
Further, Lemma~\ref{lem:dxy} guarantees that $d_\varsigma(x,y)\neq0$ for all $|x|\geq r_1$, $y\in\partial B_{r_0}$.

\begin{lemma}\label{lem:TildeUofU}
Let Assumptions~\ref{ass:TildeAlpha} and \ref{ass:r1r0} be satisfied.
Let $(\omega,u)$ be a solution to~\eqref{eq:ResonanceProblemFormal} and $\Re(i\omega d_0)<0$.
Then
\begin{align}\label{eq:DefuTilde}
\tilde u(x):=\left\{\begin{tabular}{ll}
$u(x)$, & $|x|\leq r_1$,\\
$\frac{i\omega}{4\pi} \int_{\partial B_{r_0}}
u(y) \nabla_y \tilde G_{\omega,\varsigma}(x,y)\cdot\nv(y) - \tilde G_{\omega,\varsigma}(x,y) \nabla_y u(y)\cdot\nv(y) \dd y$, &$|x|>r_1$,
\end{tabular}\right.
\end{align}
with $\tilde G_{\omega,\varsigma}(x,y):=\det\varsigma^{-1/2} h^{(1)}_0(\omega d_\varsigma(x,y))$
is well defined and $\tilde u\in H^1_0(\Omega)$.
In addition there exist constants $c_1, c_2>0$ such that $|\tilde u(x)|, |\nabla \tilde u(x)| \leq c_1 e^{-c_2|x|}$ for all $|x|>r_1$.
Let $d_\infty:=\lim_{r\to\infty} \tilde d(r)$.
For each $\epsilon>0$ the constant $c_2$ be chosen such that $c_2>-\Re(i\omega d_0)|d_\infty|\varsigma_\mathrm{max}^{-1}-\epsilon$.
\end{lemma}
\begin{proof}
Let $R_0, R_1>0$ be such that $r_1>R_1>r_0>R_0$ and $\Omega^c\subset B_{R_0}$.
It follows with standard regularity theory that $u\in H^k(B_{R_1}\setminus\ol{B_{R_0}})$ for all $k>0$.
Then the Sobolev embedding theorem yields that $u\in C^1(B_{R_1}\setminus\ol{B_{R_0}})$.
Thus it follows with Lemma~\ref{lem:dxy} that the integral in \eqref{eq:DefuTilde} is well-defined and $\tilde u\in \tilde H^1_\mathrm{loc}(B_{R_1}\setminus\ol{B_{R_0}})$.
Note that $\lim_{|x|\to\infty} d_\varsigma(x,y)/\tilde r_x=1$ uniformly in $y\in \partial B_{r_0}$ and $\lim_{|x|\to\infty} \tilde r_x/r_x=|d_\infty|d_0$.
Thus for each $\epsilon>0$ there exist constants $C_1>0$, $C_2>-\Re(i\omega d_0)|d_\infty|\varsigma_\mathrm{max}^{-1}-\epsilon$ such that $|e^{i\omega d_\varsigma(x,y)}| \leq C_1 e^{-C_2 |x|}$ for all $|x|>r_1$, $|y|=r_0$.
The functions $|\tilde u(x)|$ and $|\nabla\tilde u(x)|$ can be estimated in terms of $|\tilde G_{\omega,\varsigma}(x,y)|, |\nabla_x \tilde G_{\omega,\varsigma}(x,y)|, |\nabla_y \tilde G_{\omega,\varsigma}(x,y)|$ and $|\nabla_x\nabla_y \tilde G_{\omega,\varsigma}(x,y)|$, which itselves can be estimated in terms of $|e^{i\omega d_\varsigma(x,y)}|$.
Thus there exist constants $c_1, c_2>0$ such that $|\tilde u(x)|, |\nabla \tilde u(x)| \leq c_1 e^{-c_2|x|}$ for all $|x|>r_1$ and hence $\tilde u\in H^1(B_{R_1}\setminus\ol{B_{R_0}})$.
At $|x|=r_1$, $|y|=r_0$ it holds $\tilde G_{\omega,\varsigma}(x,y)=G_{\omega,\varsigma}(x,y)$ and $\nabla_y\tilde G_{\omega,\varsigma}(x,y)=\nabla_y G_{\omega,\varsigma}(x,y)$.
Thus the Dirichlet traces coincide at $\partial B_{r_0}$ and hence $\tilde u\in H^1_0(\Omega)$.
\end{proof}
Next we derive an equation for $\tilde u$. To this end let
\begin{align*}
\tilde\varsigma&:=(\tilde d^2d)^{-1}FMF^*\varsigma FMF^*,
\quad
F:=(\hat x,\partial_\phi \hat x,\partial_\theta \hat x),
\quad
M:=\diag(\tilde d^2,\tilde d d,\tilde d d).
\end{align*}
Note that $F$ is a real orthogonal matrix, i.e.\ $F^{-1}=F^*=F^\top$.
Let
\begin{align*}
a(\omega;u,u'):=\spl \tilde\varsigma \nabla u,\nabla u'\spr_{L^2(\Omega)}-\omega^2\spl \tilde d^2d u,u'\spr_{L^2(\Omega)}
\end{align*}
and consider the eigenvalue problem to
\begin{align}\label{eq:ResonanceProblemVariational}
\text{find }(\omega,\tilde u)\in\setC\times H^1_0(\Omega)\setminus\{0\}\text{ such that}\quad a(\omega;\tilde u,u')=0
\quad\text{for all }u'\in H^1_0(\Omega).
\end{align}

\begin{lemma}\label{lem:TildeUsolution}
Let Assumptions~\ref{ass:TildeAlpha} and \ref{ass:r1r0} be satisfied.
For a solution $(\omega,u)$ to~\eqref{eq:ResonanceProblemFormal} such that $\Re(i\omega d_0)<0$ let $\tilde u$ be as defined in Lemma~\ref{lem:TildeUofU}.
Then $(\omega,\tilde u)$ is a solution to~\eqref{eq:ResonanceProblemVariational}.
\end{lemma}
\begin{proof}
In order to show that $(\omega,\tilde u)$ solves~\eqref{eq:ResonanceProblemVariational} we apply the technique of~\cite[Theorem~4.1]{BramblePasciakTrenev:10}.
We make the dependency of quantities on $\gamma$ explicit with an index $\gamma$.
It can be seen as in the proof of Lemma~\ref{lem:TildeUofU} that $\gamma\mapsto\tilde u_\gamma$ is holomorphic.
Let $u'\in C_0^\infty(\Omega)$.
Then $\gamma\mapsto a_\gamma(\omega;\tilde u_\gamma,u')$ is holomorphic as well.
For $\gamma\in\setR_{\geq0}$ the map $\tilde x_\gamma(x):=\tilde r_\gamma(|x|) \hat x$ is real domain transformation $\Omega\to\Omega$.
Since $\tilde r_\gamma$ is strictly monoton $\tilde x_\gamma$ is even bijective.
The transformation rule yields
\begin{align*}
F(\gamma):=a_\gamma(\omega;\tilde u_\gamma,u')=a_0(\omega;u,u'\circ\tilde x_\gamma^{-1})=\spl -\div\varsigma\nabla u-\omega^2u,u'\circ\tilde x_\gamma^{-1}\spr_{L^2(\Omega)}.
\end{align*}
Since $(\omega,u)$ solves~\eqref{eq:ResonanceProblemFormal} it follows $F(\gamma)=\spl -\div\varsigma\nabla u-\omega^2u,u'\circ\tilde x_\gamma^{-1}\spr_{L^2(\Omega)}=0$ for all $\gamma\in\setR_{\geq0}$.
Thus $\gamma\mapsto F(\gamma)$ is a holomorphic function which vanishes on $\setR_{\geq0}$.
It follows that $F(\gamma)$ vanishes everywhere.
Hence $a(\omega;u,u')=0$ for all $u'\in C_0^\infty(\Omega)$.
Since $C_0^\infty(\Omega)$ is dense in $H^1_0(\Omega)$, it follows $a(\omega;u,u')=0$ for all $u'\in H^1_0(\Omega)$.
It remains to show $\tilde u\neq0$.
Since with $G_{\omega,\varsigma}$ there exists a Green's function for $-\div\varsigma \nabla u-\omega^2 u=0$ it follows that $u|_{B_{r_1}\cap \Omega}=0$ would imply $u=0$.
Thus $u|_{B_{r_1}\cap \Omega}\neq0$ and hence $\tilde u\neq0$.
It follows that $(\omega,\tilde u)$ solves~\eqref{eq:ResonanceProblemVariational}.
\end{proof}

Note that it follows from Lemmata~\ref{lem:TildeUofU} and \ref{lem:TildeUsolution} that each eigenvalue $\omega$ of \eqref{eq:ResonanceProblemFormal} with $\Re(i\omega d_0)<0$ is an eigenvalue of \eqref{eq:ResonanceProblemVariational}.
Vice-versa, for a solution $(\omega,\tilde u)$ to~\eqref{eq:ResonanceProblemVariational} such that $\Re(i\omega d_0)<0$ we can define
\begin{align*}
u(x):=\left\{\begin{tabular}{ll}
$\tilde u(x)$, & $|x|\leq r_1$,\\
$\frac{i\omega}{4\pi} \int_{\partial B_{r_0}}
\tilde u(y) \nabla_y G_{\omega,\varsigma}(x,y)\cdot\nv(y) - G_{\omega,\varsigma}(x,y) \nabla_y \tilde u(y)\cdot\nv(y) \dd y$, &$|x|>r_1$,
\end{tabular}\right..
\end{align*}
Then $(\omega,u)$ solves to~\eqref{eq:ResonanceProblemFormalA}-\eqref{eq:ResonanceProblemFormalC}.
However, we cannot guarantee that $d_\varsigma(x,y)\neq0$ for arbitrary $x\neq y$ and hence it is not clear that $\tilde G_{\omega,\varsigma}$ is a Green's function for $-\div\tilde\varsigma \nabla u-\omega^2 \tilde d^2d u=0$.
Also, we are not aware of a unique continuation principle for $-\div\tilde\varsigma \nabla u-\omega^2 \tilde d^2d u=0$, because common results require real valued coefficients.
Thus we cannot guarantee that $\tilde u|_{B_{r_1}\cap\Omega}\neq0$, which would imply $u\neq0$.

\section{Weak T-coercivity}\label{sec:wTc}

Let
\begin{align*}
\spl A(\omega)u,u'\spr_{H^1(\Omega)}&:=a(\omega;u,u'), \quad u,u'\in H^1_0(\Omega),\\
\Lambda_{d_0}&:=\{z\in\setC\colon \Re(izd_0)\neq0\}.
\end{align*}
Our aim in this section is to show that $A(\cdot)$ is weakly $T(\cdot)$-coercive on $\Lambda_{d_0}$.
By that we mean that for each $\omega\in\Lambda_{d_0}$ there exists a compact operator $K(\omega)\in L(H^1_0(\Omega))$ such that
$\inf_{u\in H^1_0(\Omega)\setminus\{0\}} |\spl (T(\omega)^*A(\omega)+K(\omega)u,u\spr_{H^1(\Omega)}|/\|u\|_{H^1(\Omega)}^2>0$.
The main difference in the analysis compared to the isotropic case \cite{Halla:19Diss,Halla:20PML} is Lemma~\ref{lem:NumRan} and an additional assumption, which involves the scaling parameters and the degree of anisotropicy $1-\frac{\varsigma_\mathrm{min}}{\varsigma_\mathrm{min}}$.
As in the isotropic case \cite{Halla:19Diss,Halla:20PML} we require an additional assumption on the profile function $\tilde\alpha$.

\begin{assumption}\label{ass:TildeAlpha2}
Let $\tilde\alpha$, $\gamma$ and $r_1$ be as is Assumption~\ref{ass:TildeAlpha} and $\tilde d$, $d$ be as in~\eqref{eq:ComplexScalingQuantities}. Let
\begin{enumerate}
 \item\label{it:limtdd} $\displaystyle\lim_{r\to+\infty} \tilde d(r)|d(r)|/\big(|\tilde d(r)|d(r)\big)=1$,
 \item\label{it:limdabsd} $\displaystyle\lim_{r\to+\infty} \big(\partial_r(\tilde d/|\tilde d|)\big)(r)
 =\lim_{r\to\infty} \big(\partial_r(d/|d|)\big)(r)=0$.
\end{enumerate}
\end{assumption}

In the isotropic case \cite{Halla:19Diss,Halla:20PML} the matrix function $\tilde\varsigma$ is diagonal, which allows to analyze its numerical range rather easily.
In the anisotropic case $\tilde\varsigma$ is not diagonal anymore.
We introduce the next lemma to analyze its numerical range.

\begin{lemma}\label{lem:NumRan}
Let $B\in\setR^{3\times 3}$ be symmetric positive definite, $\tau\in(-\pi/2,\pi/2)$ and
\begin{align*}
B_{\tau}:=\bpm e^{i\tau} B_{11} & B_{12} & B_{13} \\
B_{21} & e^{-i\tau} B_{22} & e^{-i\tau} B_{23} \\
B_{31} & e^{-i\tau} B_{32} & e^{-i\tau} B_{33} \epm.
\end{align*}
Let $\lambda_\mathrm{min}$ and $\lambda_\mathrm{max}$ be the minimal and maximal eigenvalues of $B$.
Then for all $x\in\setC^3$ with $|x|=1$ it holds
\begin{align*}
\Re(x^* B_{\tau} x) &\in [\lambda_\mathrm{min}-(1-\cos\tau)\lambda_\mathrm{max}, \,\, \lambda_\mathrm{max}-(1-\cos\tau)\lambda_\mathrm{min}],\\
\Im(x^* B_{\tau} x) &\in [-\lambda_\mathrm{max}|\sin\tau|,\,\, \lambda_\mathrm{max}|\sin\tau|].
\end{align*}
\end{lemma}
\begin{proof}
We employ the representation $B_\tau=B-(1-\cos\tau)\tilde B_+ +i\sin\tau\tilde B_-$ with
\begin{align*}
\tilde B_\pm:=\bpm B_{11} & 0& 0 \\
0 & \pm B_{22} & \pm B_{23} \\
0 & \pm B_{32} & \pm B_{33} \epm.
\end{align*}
Since $B, \tilde B_+$ and $\tilde B_-$ are symmetric it follows that
$\Im(x^* B_{\tau} x)=\sin\tau \, x^*\tilde B_- x$ and
$\Re(x^* B_{\tau} x)=x^*Bx -(1-\cos\tau) x^*\tilde B_+ x$.
We estimate
\begin{align*}
|x^*\tilde B_-x|
&=\left| B_1|x_1|^2 -\bpm x_2\\x_3 \epm^* \bpm B_{22}&B_{23}\\B_{32}&B_{33} \epm \bpm x_2\\x_3 \epm \right|\\
&\leq B_1|x_1|^2 +\bpm x_2\\x_3 \epm^* \bpm B_{22}&B_{23}\\B_{32}&B_{33} \epm \bpm x_2\\x_3 \epm
= x^*\tilde B_+ x.
\end{align*}
Further, we compute
\begin{align*}
&\sup_{x\in\setC^3, |x|=1} x^*\tilde B_+ x \\
&= \sup_{x\in\setC^3, |x|=1}
B_1|x_1|^2 +\bpm x_2\\x_3 \epm^* \bpm B_{22}&B_{23}\\B_{32}&B_{33} \epm \bpm x_2\\x_3 \epm \\
&\leq \sup_{x\in\setC^3, |x|=1}
B_1|x_1|^2 +\bpm x_2\\x_3 \epm^* \bpm B_{22}&B_{23}\\B_{32}&B_{33} \epm \bpm x_2\\x_3 \epm
+2\big|\ol{x_1}(B_{12}x_2+B_{13}x_3)\big| \\
&= \sup_{x\in\setC^3, |x|=1}
B_1|x_1|^2 +\bpm x_2\\x_3 \epm^* \bpm B_{22}&B_{23}\\B_{32}&B_{33} \epm \bpm x_2\\x_3 \epm
+2\Re\big( \ol{x_1}(B_{12}x_2+B_{13}x_3) \big) \\
&=\sup_{x\in\setC^3, |x|=1} x^*Bx
\leq \lambda_\mathrm{max}.
\end{align*}
Like-wise we obtain $\inf_{x\in\setC^3, |x|=1} x^*\tilde B_+ x\geq \lambda_\mathrm{min}$.
So $|\Im(x^* B_{\tau} x)| \leq \lambda_\mathrm{max}|\sin\tau|$.
In addition, we estimate for $x\in\setC^3$ with $|x|=1$ that
\begin{align*}
x^*Bx -(1-\cos\tau) x^*\tilde B_+ x \geq \lambda_\mathrm{min}-(1-\cos\tau)\lambda_\mathrm{max},\\
x^*Bx -(1-\cos\tau) x^*\tilde B_+ x \leq \lambda_\mathrm{max}-(1-\cos\tau)\lambda_\mathrm{min},
\end{align*}
which proves the claim.
\end{proof}

We note that $\lambda_\mathrm{min}-(1-\cos\tau)\lambda_\mathrm{max}>0$ if and only if
$\cos\tau>1-\lambda_\mathrm{min}/\lambda_\mathrm{max}$.
Henceforth, we consider the logarithmic function $\log$ with the branch cut $\setR_{<0}$, i.e.\ we define the argument of a complex number as
\begin{align*}
\arg z \colon \setC\setminus\{0\}\to [-\pi,\pi), \quad z=|z|\exp(i\arg z).
\end{align*}
In order to formulate our analysis we introduce
\begin{align*}
H&:=\diag(|\tilde d/d|^{1/2},\, |\tilde d/d|^{-1/2},\, |\tilde d/d|^{-1/2}),\\
\tau&:=\arg (\tilde d/d),\\
\tau_*&:=\sup_{r>r_1} |\arg (\tilde d/d)|,\\
\psi&:=
\arg\big(\varsigma_\mathrm{min}-(1-\cos\tau_*)\varsigma_\mathrm{max}-i\varsigma_\mathrm{max}\sin\tau \big),\\
\psi_*&:=\sup_{r>r_1} \psi.
\end{align*}
We adopt the overloaded notation \eqref{eq:OverloadedNotation} also for $H, \tau, \psi$.
The definition of $\psi$ and Lemma~\ref{lem:NumRan} yield
\begin{align*}
\psi &\geq \max_{y\in\setC^3,|y|=1} |\arg y^* (F^\top \varsigma F)_\tau y|.
\end{align*}
With the former quantities we obtain the representation
\begin{align*}
\tilde\varsigma=\tilde d F H (F^* \varsigma F)_\tau H F^*.
\end{align*}
The assumption $\cos\tau_*>1-\varsigma_\mathrm{min}/\varsigma_\mathrm{max}$ guarantees that $\psi_* <\pi/2$.
Note that for an isotropic material the matrix $\varsigma$ is a scalar multiple of the identity matrix and hence $\varsigma_\mathrm{min}=\varsigma_\mathrm{max}$.
Thus in this case the condition $\cos\tau_*>1-\varsigma_\mathrm{min}/\varsigma_\mathrm{max}=0$ reduces to $\tau_*<\pi/2$.
Further, it follows from Assumption~\ref{ass:TildeAlpha2} that $\lim_{r\to\infty} \psi(r)=0$ .
In order to construct a suitable $T$-operator we introduce
\begin{align*}
\hat\alpha(r)&:=\left\{\begin{array}{rl}\lim_{\rho\to r_1+}\alpha(\rho)&\text{for }0\leq r\leq r_1,\\
\alpha(r)&\text{for }r>r_1,\end{array}\right.,\\
\hat d(r)&:=1+i\hat\alpha(r), \quad r\geq0,\\
\hat \tau(r)&:=\arg(\tilde d(r)/\hat d(r)), \quad r\geq0,\\
\hat\psi(r)&:=
\arg\big(\varsigma_\mathrm{min}-(1-\cos\tau_*)\varsigma_\mathrm{max}-i\varsigma_\mathrm{max}\sin\hat\tau \big).
\end{align*}
We adopt the overloaded notation \eqref{eq:OverloadedNotation} also for $\hat\alpha,\hat d,\hat\tau,\hat\psi$.
Then we define the multiplication operator
\begin{align*}
T(\omega)u:=\left\{\begin{array}{ll}
\frac{|\tilde d|}{\ol{\tilde d}} e^{i\hat\psi} u &\text{for } \arg(-\omega^2d_0^2)\in [-\pi,0], \vspace{3mm}\\
\frac{|\tilde d|}{\ol{\tilde d}} e^{-i\hat\psi} u &\text{for }\arg (-\omega^2d_0)\in (0,\pi).
\end{array}\right.
\end{align*}
The next Lemma~\ref{lem:limetazero} is necessary for Lemma~\ref{lem:Tbij} and Theorem~\ref{thm:wTc}.
\begin{lemma}\label{lem:limetazero}
Let Assumptions~\ref{ass:TildeAlpha}, \ref{ass:TildeAlpha2} and $\cos\tau_*>1-\varsigma_\mathrm{min}/\varsigma_\mathrm{max}$ be satisfied. Then
\begin{align*}
\lim_{r\to\infty}\partial_r
\bigg(\frac{|\tilde d|}{\ol{\tilde d}} e^{\pm i\hat\psi}\bigg)(r)=0.
\end{align*}
\end{lemma}
\begin{proof}
The product rule yields $\partial_r
\bigg(\frac{|\tilde d|}{\ol{\tilde d}} e^{\pm i\hat\psi}\bigg)
=\partial_r\Big(\frac{|\tilde d|}{\ol{\tilde d}}\Big)e^{\pm i\hat\psi}
+\frac{|\tilde d|}{\ol{\tilde d}} \partial_r e^{\pm i\hat\psi}$.
Due to Assumption~\ref{ass:TildeAlpha2} the claim follows, if we show $\lim_{r\to\infty}\partial_r e^{\pm i\hat\psi}=0$.
It suffices to consider $r>r_1$ and thus $\hat\psi=\psi$.
We compute $\partial_r e^{\pm i\psi}=\pm ie^{\pm i\psi} \partial_r \psi$ and
$\partial_r \psi=\arctan'\Big( \frac{\varsigma_\mathrm{max}\sin\tau}{\varsigma_\mathrm{min}-(1-\cos\tau_*)\varsigma_\mathrm{max}} \Big)
\frac{\varsigma_\mathrm{max}\cos\tau}{\varsigma_\mathrm{min}-(1-\cos\tau_*)\varsigma_\mathrm{max}} \partial_r \tau$.
It follows from Assumptions~\ref{ass:TildeAlpha} and \ref{ass:TildeAlpha2} that $|\tau|\leq\tau_*<\pi/2$.
Hence it suffices to show $\lim_{r\to\infty}\partial_r \tau=0$.
Due to $\tau=\arg(\tilde d/d)=\arg\big(\frac{\tilde d}{|\tilde d|}\frac{|d|}{d}\big)$ it follows
\begin{align*}
\partial_r \tau=\arg'(\tilde d/d) \partial_r \Big(\frac{\tilde d}{|\tilde d|}\frac{|d|}{d}\Big).
\end{align*}
Since $\arg'(\tilde d/d)$ is bounded, the claim follows with the product rule and Assumption~\ref{ass:TildeAlpha2}.
\end{proof}

\begin{lemma}\label{lem:Tbij}
Let Assumptions~\ref{ass:TildeAlpha}, \ref{ass:TildeAlpha2} and $\cos\tau_*>1-\varsigma_\mathrm{min}/\varsigma_\mathrm{max}$ be satisfied.
Then $T(\omega)$ is in $L(H^1_0(\Omega))$ and bijective.
\end{lemma}
\begin{proof}
A multiplication operator with symbol $m$ such that $|m|=1$ everywhere and $\nabla m\in L^\infty(\Omega)$ is in $L(H^1_0(\Omega))$ and bijective.
Denote the symbol of the multiplication operator $T(\omega)$ as $m$.
It holds $|m|=1$ everywhere and $m$ admits a weak derivative $\nabla m$ due to Assumption~\ref{ass:TildeAlpha}.
Further, it follows from Assumption~\ref{ass:TildeAlpha} and Lemma~\ref{lem:limetazero} that $\nabla m\in L^\infty(\Omega)$.
Thus the claim is proven.
\end{proof}
The next Lemma~\ref{lem:WeigthedL2Compactness} is necessary for Theorem~\ref{thm:wTc}.
\begin{lemma}\label{lem:WeigthedL2Compactness}
Let $\eta\in L^\infty(\Omega)$ be so that $\lim_{r\to\infty}\|\eta\|_{L^\infty(\Omega\cap B_r^c)}=0$.
Then the multiplication and embedding operator $K_{\eta}\colon H^1_0(\Omega)\to L^2(\Omega)\colon u\mapsto \eta u$ is compact.
\end{lemma}
\begin{proof}
Follows from \cite[Lemma~4.3]{Halla:19Diss} or \cite[Lemma~3.3]{Halla:20PML}.
\end{proof}

\begin{theorem}\label{thm:wTc}
Let Assumptions~\ref{ass:TildeAlpha}, \ref{ass:TildeAlpha2} and $\cos\tau_*>1-\varsigma_\mathrm{min}/\varsigma_\mathrm{max}$ be satisfied.
Let $\omega\in \Lambda_{d_0}$.
Then $A(\omega)$ is weakly $T(\omega)$-coercive.
\end{theorem}
\begin{proof}
First note that $T(\omega)\in L(H^1_0(\Omega))$ is bijective due to Lemma~\ref{lem:Tbij}.
Let $\arg(-\omega^2d_0^2)\in (-\pi,0]$.
We split $T(\omega)^*A(\omega)=A_1+A_2$ into a coercive operator $A_1$ and compact operator $A_2$.
Let $\eta:=\nabla \frac{e^{i\hat\psi}|\tilde d|}{\ol{\tilde d}}$ and $A_1,A_2\in L(H^1_0(\Omega))$ be defined by
\begin{align*}
\spl A_1u,u' \spr_{H^1(\Omega)}&:=
\spl \frac{e^{-i\hat\psi}|\tilde d|}{\tilde d} \tilde\varsigma \nabla u,\nabla u' \spr_{L^2(\Omega)}
-\omega^2d_0^2 \spl |\tilde d^2d| u,u'\spr_{L^2(\Omega)},\\
\spl A_2u,u' \spr_{H^1(\Omega)}&:=
\spl \tilde\varsigma \nabla u,\eta u' \spr_{L^2(\Omega)}
-\omega^2 \spl (e^{-i\hat\psi} \tilde dd/|\tilde dd|-d_0^2) |\tilde d^2d| u,u'\spr_{L^2(\Omega)}
\end{align*}
for all $u,u'\in H^1_0(\Omega)$.
Then indeed $T(\omega)^*A(\omega)=A_1+A_2$.
Recall
\begin{align*}
\tilde\varsigma=\tilde d F H (F^* \varsigma F)_\tau H F^*.
\end{align*}
Hence
\begin{align*}
\frac{e^{-i\hat\psi}|\tilde d|}{\tilde d} \tilde\varsigma=
|\tilde d| F H e^{-i\hat\psi}(F^* \varsigma F)_\tau H F^*.
\end{align*}
For $r\geq r_1$ the numerical range of $e^{-i\hat\psi}(F^* \varsigma F)_\tau$ is contained in the salient sector spanned by $1$ and $e^{-2i\hat\psi_*}$ due to Lemma~\ref{lem:NumRan} and the definition of $\hat\psi$.
For $r<r_1$ it holds $e^{-i\hat\psi}(F^* \varsigma F)_\tau=e^{-i\hat\psi(r_1)}F^* \varsigma F$.
Thus the numerical range of $e^{-i\hat\psi}(F^* \varsigma F)_\tau$ is contained in the salient sector spanned by $1$ and $e^{-2i\hat\psi_*}$ for all $r>0$.
Let
\begin{align*}
\zeta:=\min(-2\psi_*, \arg(-\omega^2d_0^2))
\end{align*}
It follows $\zeta\in(-\pi,0)$.
Hence we can estimate
\begin{align*}
\Re(ie^{-i(\pi+\zeta)/2} \spl A_1u,u\spr_{H^1(\Omega)}) \geq C \|u\|_{H^1_0(\Omega)}^2
\end{align*}
with
\begin{align*}
C&:=\cos(\zeta/2) \min\{ C_H(\varsigma_\mathrm{min}-(1-\cos\tau_*)\varsigma_\mathrm{max}), \,\, |\omega^2|\},\\
C_H&:=\inf_{r>0} \min\{|\tilde d|, |\tilde d^2/d|, |d|\},
\end{align*}
whereby we used that $|\tilde d|,|d|\geq1$.
Note that by assumption of this theorem $\varsigma_\mathrm{min}-(1-\cos\tau_*)\varsigma_\mathrm{max}>0$.
Further, $C_H>0$ due to Assumption~\ref{ass:TildeAlpha} and $\omega\neq0$ due to $\omega\in\Lambda_{d_0}$.
Thus $C>0$ and hence $A_1$ is coercive.

Let us now analyze $A_2$.
Due to Lemma~\ref{lem:limetazero} it holds $\lim_{r\to\infty} \|\eta\|_{L^\infty(B_r^c)}=0$.
Further, $\lim_{r\to\infty} (e^{-i\hat\psi}\tilde dd/|\tilde dd|-d_0^2) |\tilde d^2d|=0$ due to Assumptions~\ref{ass:TildeAlpha}, \ref{ass:TildeAlpha2} and the definition of $d_0$.
Hence, it follows with Lemma~\ref{lem:WeigthedL2Compactness} that $A_2$ is compact.
Thus $A(\omega)$ is weakly $T(\omega)$-coercive for $\arg(-\omega^2d_0^2)\in (-\pi,0]$.
The proof for $\arg(-\omega^2d_0^2)\in (0,\pi)$ can be seen similarly as for $\arg(-\omega^2d_0^2)\in [-\pi,0]$.
\end{proof}

Assumptions~\ref{ass:TildeAlpha} and \ref{ass:TildeAlpha2} already appear for isotropic materials.
They can be checked easily for common profiles.
We emphasize that Assumption~\ref{ass:r1r0} is not required for Theorem~\ref{thm:wTc}.
On the other hand, the assumption $\cos\tau_*>1-\varsigma_\mathrm{min}/\varsigma_\mathrm{max}$ deserves some discussion.
Recall $\tau_*=\sup_{r>r_1}|\tau|$, $|\tau|=\arg(d/\tilde d)$ and
\begin{align*}
d/\tilde d
=\frac{1+\gamma\tilde\alpha+\gamma r\partial_r\tilde\alpha}{1+\gamma\tilde\alpha}
=1+\frac{\gamma r\partial_r\tilde\alpha}{1+\gamma\tilde\alpha}
=1+\frac{1}{1/\gamma+\frac{\tilde\alpha}{r\partial_r\tilde\alpha}}
\end{align*}
Hence the argument of $d/\tilde d$ decreases as $\Im(\gamma)$ decreases or $\frac{\tilde\alpha}{r\partial_r\tilde\alpha}$ increases.
For the affine profile
\begin{align*}
\tilde\alpha(r)=1-r_1/r, \qquad
\tilde r(r)=\gamma(r-r_1)+r
\end{align*}
it holds $d=1+\gamma$ and hence $\tau_*$ can be computed explicitly as $\tau_*=\arg(1+\gamma)$.
This means that the measure of anisotropicy $1-\varsigma_\mathrm{min}/\varsigma_\mathrm{max}$ limits $\arg(1+\gamma)$.
For fixed $\Re(\gamma)$ this limits $\Im(\gamma)$, i.e.\ the rate of the damping.

To derive equations for the transient wave equation, the constant $\gamma$ is chosen $\omega$-dependent.
The most convenient choice is $\gamma(\omega)=\frac{1}{-i\omega}$.
However, for an affine scaling this leads to $\tau_*=\arg\big(1+\frac{1}{-i\omega}\big)$ and thus $\lim_{\omega\to0}\tau_*=+\infty$.
Hence, we do expect that the derived system admits instabilities as the slow waves reach the complex scaled layer.
As a remedy we propose to choose $\gamma(\omega)=\frac{1}{c-i\omega}$ with a positive constant $c$.
This way we obtain $\tau_*=\arg\big(1+\frac{1}{c-i\omega}\big)\leq\arctan\Big(\frac{1}{2\sqrt{c^2+c}}\bigg)$.
Hence $\tau_*$ can be bounded uniformly for all frequencies $\omega\in\setR$ and be made arbitrarily small as $c\to+\infty$.
Thus for an anisotropic material with matrix $\varsigma$ we choose $c$ big enough such that
\begin{align*}
\cos\arctan\Big(\frac{1}{2\sqrt{c^2+c}}\bigg)=\frac{1}{\sqrt{1+\frac{1}{4(c^2+c)}}}>1-\varsigma_\mathrm{min}/\varsigma_\mathrm{max}.
\end{align*}
Note that with this scaling the damping of waves $\Re(\gamma(\omega)i\omega)=\frac{-\omega^2}{c^2+\omega^2}$ is not uniform in the frequency and vanishes for $\omega\to0$.
However, this seems unavoidable.

\section{Convergence of approximations}\label{sec:Appr}

To analyze the approximation of \eqref{eq:ResonanceProblemVariational} we can proceed exactly as in \cite{Halla:19Diss,Halla:20PML}, and hence we keep our presentation brief.
This correspondence of the analysis' of approximations underlines our argument that for a radial CS/PML the difference between isotropic and anisotropic equations is marginal.

The convenient way to approximate \eqref{eq:ResonanceProblemVariational} consists of two steps.
First, one chooses a sequence of subdomains $(\Omega_n)_{n\in\setN}$ such that for each $n\in\setN$
\begin{enumerate}
 \item $\Omega_n$ is a bounded Lipschitz domain,
 \item $\Omega_n\subset\Omega$,
 \item $\partial\Omega\subset\partial\Omega_n$,
 \item $\partial\Omega_n\setminus\partial\Omega$ splits $\Omega$ into two connected parts,
 \item for each $R>0$ exists an index $n_0\in\setN$ such that $(\Omega\cap B_R)\subset \Omega_n$
for all $n>n_0$.
\end{enumerate}
Subsequently one imposes a homogeneous Dirichlet boundary condition on $\partial\Omega_n\setminus\partial\Omega$ and considers the problem on the truncated domain $\Omega_n$:
\begin{align}\label{eq:ResonanceProblemTruncated}
\begin{aligned}
\text{find }(\omega,u)\in\setC\times H^1_0(\Omega_n)\setminus\{0\}&\text{ such that}\\
\spl \tilde\varsigma \nabla u,\nabla u'\spr_{L^2(\Omega_n)}-\omega^2\spl &\tilde d^2d u,u'\spr_{L^2(\Omega_n)}=0
\quad\text{for all }u'\in H^1_0(\Omega_n).
\end{aligned}
\end{align}
In a second step one chooses a fixed $n\in\setN$ and discretizes problem~\eqref{eq:ResonanceProblemTruncated} by a Galerkin approximation with finite element spaces $X_h(\Omega_n)\subset H^1_0(\Omega_n)$.
In order to guarantee that for any combination of $n\to\infty,h\to0+$ the approximations converge we follow the approach of \cite{HohageNannen:15,Halla:19Diss,Halla:20PML}.
That is we identify a function $u\in H^1_0(\Omega_n)$ with its extension by zero in $\Omega\setminus\Omega_n$: $\hat u|_{\Omega_n}=u$ and $\hat u|_{\Omega\setminus\Omega_n}=0$.
Obviously it follows that $\hat u\in H^1_0(\Omega)$.
Then the truncated discretized problem can be formulated as
\begin{align*}
\text{find }(\omega,u)\in\setC\times \hat X_h(\Omega_n)\setminus\{0\}&\text{ such that}\quad
a(\omega;u,u')=0
\quad\text{for all }u'\in \hat X_h(\Omega_n)
\end{align*}
with $\hat X_h(\Omega_n):=\{\hat u\colon u\in X_h(\Omega_n)\}\subset H^1_0(\Omega)$.
Thus we can analyze the approximations through simultaneous domain truncation and discretization inside the convenient framework of Galerkin approximations.
Hence, from now on we consider a sequence of finite dimensional subspaces $(X_n\subset H^1_0(\Omega))_{n\in\setN}$ with associated orthogonal projections $P_n$ such that $\lim_{n\to\infty}\|u-P_nu\|_{H^1(\Omega)}=0$ for each $u\in H^1_0(\Omega)$.
We say that the Galerkin approximation $A_n(\cdot):=P_nA(\cdot)|_{X_n}$ of $A(\cdot)\colon\Lambda_{d_0}\to L(H^1_0(\Omega))$ is $T(\cdot)$-compatible, if there exist operator functions $T_n(\cdot)\colon\Lambda_{d_0}\to L(X_n)$ such that $\lim_{n\to\infty} \|T(\omega)-T_n(\omega)\|_n=0$ for each $\omega\in\Lambda_{d_0}$, whereby
\begin{align*}
\|T(\omega)-T_n(\omega)\|_n:=
\sup_{u_n\in X_n\setminus\{0\}}\frac{\|(T(\omega)-T_n(\omega))u_n\|_{H^1(\Omega)}}{\|u_n\|_{H^1(\Omega)}}.
\end{align*}
For such approximations we can employ the framework of $T$-compatible approximations of weakly $T$-coercive operators \cite{Halla:19Tcomp,Halla:19Diss} which yields the convergence of eigenvalues and eigenfunctions.
It remains to construct the operators $T_n(\omega)$.
To this end we first introduce a slight modification of $T(\cdot)$.
For $\eta\in W^{1,\infty}(\Omega)$ we denote $M_\eta\in L(H^1_0(\Omega))$ the multiplication operator with symbol $\eta$.
Let $\eta$ be the symbol of the multiplication operator $T(\omega)$.
Then it follows from \cite[Lemma~4.13]{Halla:19Diss} that for each $\epsilon>0$ exist $\eta_\epsilon\colon (0,+\infty)\to\setC$ and $\hat r_1, \hat r_2\in (r_1,+\infty)$ such that
\begin{enumerate}
 \item $\|\eta-\eta_\epsilon\|_{L^\infty(0,+\infty)}<\epsilon$,
 \item $\eta_\epsilon(r)=\eta(r_1)$ for $r \leq \hat r_1$,
 \item $\eta_\epsilon(r)=d_0$ for $r \geq \hat r_2$,
 \item $\eta_\epsilon$ is infinitely many times differentiable.
\end{enumerate}
A slight adaptation of the proof of Theorem~\ref{thm:wTc} yields that for each $\omega\in\Lambda_{d_0}$ we can find $\epsilon(\omega)>0$ such that $A(\omega)$ is weakly $\tilde T(\omega):=M_{\eta_{\epsilon(\omega)}}$-coercive.
For our analysis we require two additional assumptions on the Galerkin spaces $X_n$.

\begin{assumption}\label{ass:LocalApproximation}
There exists a sequence $\big(h(n)\big)_{n\in\setN}\in (\setR^+)^\setN$ with $\lim_{n\in\setN} h(n)$ $=$ $0$.
There exist bounded linear projection operators $\Pi_n\colon H^1_0(\Omega)\to X_n, n\in\setN$
that act locally in the following sense: there exist constants $C_1, R^*>1$ such that for $n\in\setN$, $s \in \{1,2\}$,
$x_0\in\Omega$, if $B_{R^*h(n)}(x_0)\subset \Omega$, $u\in H$ and $u|_{B_{R^*h(n)}(x_0)} \in H^s(B_{R^*h(n)}(x_0))$, then
\begin{align*}
\|u-\Pi_n u\|_{H^1(B_{h(n)}(x_0))} \leq C_1 h(n)^{s-1} \|u\|_{H^s(B_{R^*h(n)}(x_0))}.
\end{align*}
\end{assumption}

\begin{assumption}\label{ass:ConstantFunctions}
For any $D\subset\Omega$ which is compact in $\Omega$, there exists $n_0>0$
such that for each $n\in\setN, n>n_0$ there exists $u_{D,n}\in X_n$ with $u_{D,n}|_D=1$.
\end{assumption}
Both assumptions are satisfied by common finite element spaces:
for the first Assumption~\ref{ass:LocalApproximation} we can employ the Scott-Zhang interpolant (see, e.g.\ \cite[Lemma~1.130]{ErnGuermond:04}) and the second Assumption~\ref{ass:ConstantFunctions} is satisfied naturally.
With $\Pi_n$ as in Assumption~\ref{ass:LocalApproximation} we define $T_n(\omega):=\Pi_n\tilde T(\omega)|_{X_n}$.
Then it follows from \cite[Theorem~4.17]{Halla:19Diss} that $\lim_{n\to\infty} \|\tilde T(\omega)-T_n(\omega)\|_n=0$ for each $\omega\in\Lambda_{d_0}$, and hence the Galerkin approximation $A_n(\cdot)$, $n\in\setN$ is $\tilde T(\cdot)$-compatible.
We conclude in the following theorem.

\begin{theorem}\label{thm:conv}
Let Assumptions~\ref{ass:TildeAlpha}, \ref{ass:TildeAlpha2}, \ref{ass:LocalApproximation}, \ref{ass:ConstantFunctions} and $\cos\tau_*>1-\varsigma_\mathrm{min}/\varsigma_\mathrm{max}$ be satisfied, and assume that the resolvent set of $A(\cdot)\colon\Lambda_{d_0}\to L(H^1_0(\Omega))$ is not empty.
Then the eigenvalues and eigenfunctions of $A_n(\cdot)$ converge to those of $A(\cdot)$ in the following sense.
\begin{enumerate}[i)]
 \item\label{item:SP-a} For every eigenvalue $\omega_0$ of $A(\cdot)$ there exists a sequence $(\omega_n)_{n\in\setN}$
 converging to $\omega_0$ with $\omega_n$ being an eigenvalue of $A_n(\cdot)$ for almost all $n\in\setN$.
 \item\label{item:SP-b} \sloppy Let $(\omega_n, u_n)_{n\in\setN}$ be a sequence of normalized eigenpairs of $A_n(\cdot)$, i.e.\ $A_n(\omega_n)u_n$ $=0$ and $\|u_n\|_X=1$, so that $\omega_n\to \omega_0\in\Lambda_{d_0}$. Then $\omega_0$ is an eigenvalue of $A(\cdot)$, and $(u_n)_{n\in\setN}$ is a compact sequence and its cluster points are normalized eigenelements of $A(\omega_0)$.
 \item\label{item:SP-c} For every compact $\tilde\Lambda\subset\rho(A(\cdot))$ the sequence $(A_n(\cdot))_{n\in\setN}$ is  stable on $\tilde\Lambda$, i.e.\ there exist $n_0\in\setN$ and $c>0$ such that $\|A_n(\omega)^{-1}\|_{L(X_n)}\leq c$ for  all $n>n_0$ and all $\omega\in\tilde\Lambda$.
 \item\label{item:SP-d} \label{item:Stability} For every compact $\tilde\Lambda\subset\Lambda_{d_0}$ with $\tilde\Lambda\cap\sigma\big(A(\cdot)\big)=\{\omega_0\}$ and rectifiable boundary  $\partial\tilde\Lambda\subset\rho\big(A(\cdot)\big)$ there exists an index $n_0\in\setN$ such that
 \begin{align*}
 \dim G(A(\cdot),\omega_0) = \sum_{\omega_n\in\sigma\left(A_n(\cdot)\right)\cap\tilde\Lambda}
 \dim G(A_n(\cdot),\omega_n).
 \end{align*}
 for all $n>n_0$, whereby $G(B(\cdot),\omega)$ denotes the generalized eigenspace of an operator function $B(\cdot)$ at $\omega\in\Lambda_{d_0}$.
\end{enumerate}
Let $\tilde\Lambda\subset\Lambda_{d_0}$ be a compact set with rectifiable boundary $\partial\tilde\Lambda\subset\rho\big(A(\cdot)\big)$, $\tilde\Lambda\cap\sigma\big(A(\cdot)\big)=\{\omega_0\}$ and
\begin{align*}
\delta_n&:=\max_{\substack{u_0\in G(A(\cdot),\omega_0)\\\|u_0\|_{H^1(\Omega)}\leq1}} \, \inf_{u_n\in X_n} \|u_0-u_n\|_{H^1(\Omega)},\\
\delta_n^*&:=\max_{\substack{u_0\in G(A^*(\ol{\cdot}),\omega_0)\\\|u_0\|_{H^1(\Omega)}\leq1}} \, \inf_{u_n\in X_n} \|u_0-u_n\|_{H^1(\Omega)},
\end{align*}
whereby $\overline{\omega_0}$ denotes the complex conjugate of $\omega_0$ and $A^*(\cdot)$
the adjoint operator function of $A(\cdot)$ defined by $A^*(\omega):=A(\omega)^*$ for each $\omega\in\Lambda_{d_0}$.
Then there exist $n\in\setN$ and $c>0$ such that for all $n>n_0$
\begin{enumerate}[i)]
\setcounter{enumi}{4}
  \item\label{item:SP-e}
  \begin{align*}
  |\omega_0-\omega_n|\leq c(\delta_n\delta_n^*)^{1/\varkappa\left(A(\cdot),\omega_0\right)}
  \end{align*}
  for all $\omega_n\in\sigma\big(A_n(\cdot)\big)\cap\tilde\Lambda$, whereby
  $\varkappa\left(A(\cdot),\omega_0\right)$ denotes the maximal length of a Jordan
  chain of $A(\cdot)$ at the eigenvalue $\omega_0$,
  \item\label{item:SP-f}
  \begin{align*}
  |\omega_0-\omega_n^\mathrm{mean}|\leq c\delta_n\delta_n^*
  \end{align*}
  whereby $\omega_n^\mathrm{mean}$ is the weighted mean of all the eigenvalues of $A_n(\cdot)$ in $\tilde\Lambda$
  \begin{align*}
  \omega_n^\mathrm{mean}:=\sum_{\omega\in\sigma\left(A_n(\cdot)\right)\cap\tilde\Lambda_{d_0}}\omega\,
  \frac{\dim G(A_n(\cdot),\omega)}{\dim G(A(\cdot),\omega_0)},
  \end{align*}
  \item\label{item:SP-g}
  \begin{align*}
  \begin{split}
  \inf_{u_0\in\ker A(\omega_0)} \|u_n-u_0\|_X &\leq c \Big(|\omega_n-\omega_0|+
  \max_{\substack{u'_0\in\ker A(\omega_0)\\\|u_0'\|_X\leq1}} \inf_{u'_n\in X_n} \|u'_0-u'_n\|_X\Big)\\
  &\leq c\big(c(\delta_n\delta_n^*)^{1/\varkappa\left(A(\cdot),\omega_0\right)} + \delta_n\big)
  \end{split}
  \end{align*}
  for all $\omega_n\in\sigma\big(A_n(\cdot)\big)\cap\tilde\Lambda$ and all $u_n\in \ker A_n(\omega_n)$ with $\|u_n\|_X=1$.
\end{enumerate}
\end{theorem}
\begin{proof}
$A(\omega)$ is an affine function in $\omega^2$ and hence holomorphic.
As argued in this section $A(\cdot)$ is weakly $\tilde T(\cdot)$-coercive.
For each $n\in\setN$ the space $X_n$ is finite dimensional and hence the operator $A_n(\omega)$ is Fredholm with index zero.
As argued in this section there exist operators $T_n(\omega)\in L(X_n)$ with $\lim_{n\to\infty} \|\tilde T(\omega)-T_n(\omega)\|_n=0$.
Hence \cite[Theorem~3.17]{Halla:19Diss} yields the claim.
\end{proof}
Theorem~\ref{thm:conv} yields convergence rates in terms of the best approximation error $\inf_{u'\in X_n} \|u-u'\|_{H^1(\Omega)}$ with eigenfunction $u\in H^1_0(\Omega)$.
With a simple technique we can estimate this abstract term as ``truncation error+discretization error''.
To this end let $\chi\colon\setR\to[0,1]$ be infinitely many times differntiable, $\chi(r)=1$ for $r>1$ and $\chi(r)=0$ for $r<0$.
Let $B_{R+1}\cap\Omega\subset\Omega_n$ and $\hat u(x):=\chi(1+R-|x|)u(x)$.
Then we can obtain with Lemma~\ref{lem:TildeUofU} that
\begin{align*}
\|u-\hat u\|_{H^1(\Omega)}\leq (1+\|\chi\|_{W^{1,\infty}(\Omega)}) \|u\|_{H^1(B_R^c)}
\lesssim e^{-c_2R}
\end{align*}
and hence
\begin{align*}
\inf_{u'\in X_n} \|u-u'\|_{H^1(\Omega)} \lesssim e^{-c_2R} + \inf_{u'\in X_n} \|\hat u-u'\|_{H^1(\Omega)}.
\end{align*}
The term $\inf_{u'\in X_n} \|\hat u-u'\|_{H^1(\Omega)}$ is now the discretization error which can for $p$-order finite element methods be conveniently estimated as $h^p$.

\section{Computational experiments}\label{sec:comp}

\begin{figure}
\begin{subfigure}{.47\textwidth}
\centering
\includegraphics[clip, trim={20cm 2cm 20cm 2cm}, scale=0.18]{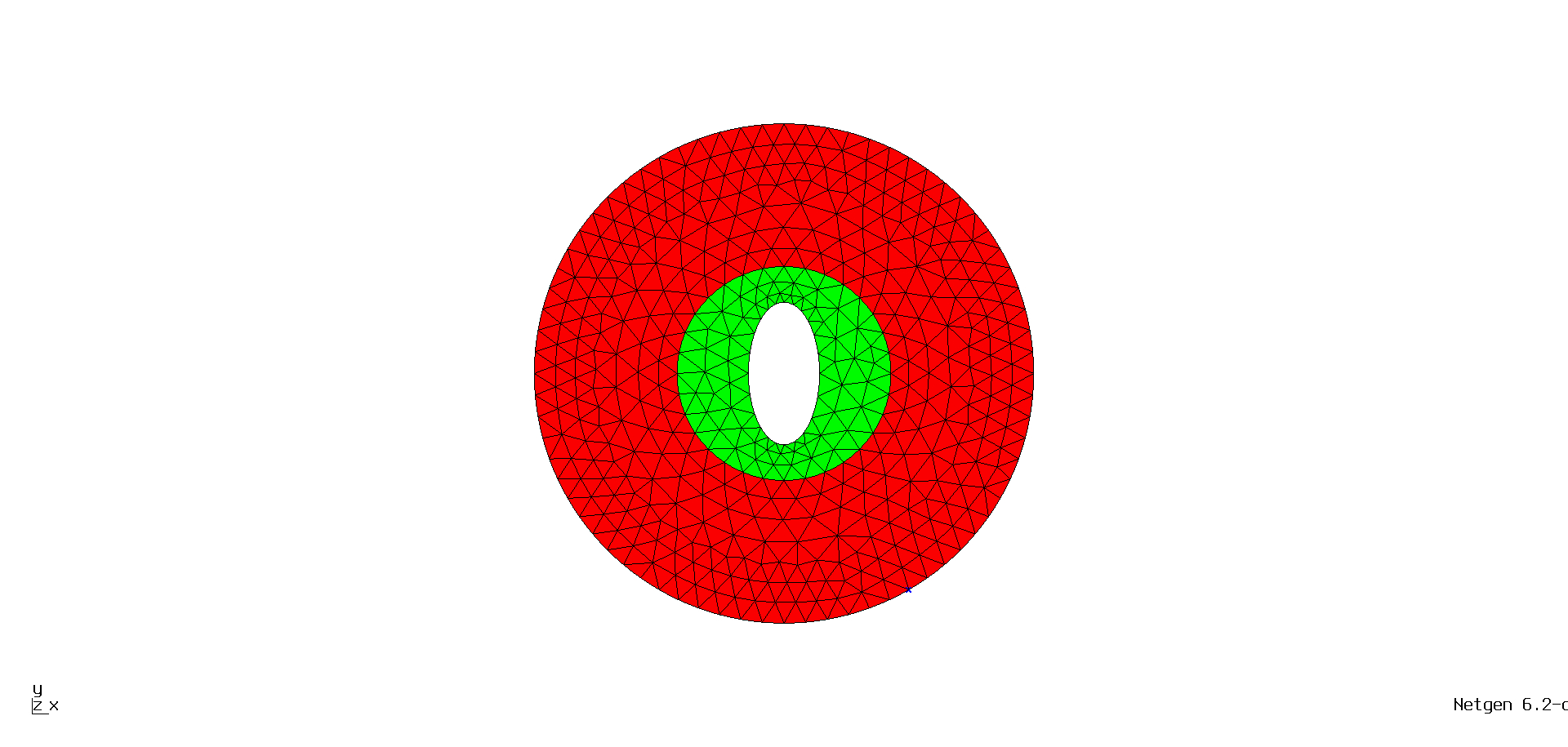}
\caption{Anisotropic case.}
\label{fig:mesh_iso}
\end{subfigure}
\hfill
\begin{subfigure}{.47\textwidth}
\centering
\includegraphics[clip, trim={20cm 2cm 20cm 2cm}, scale=0.18]{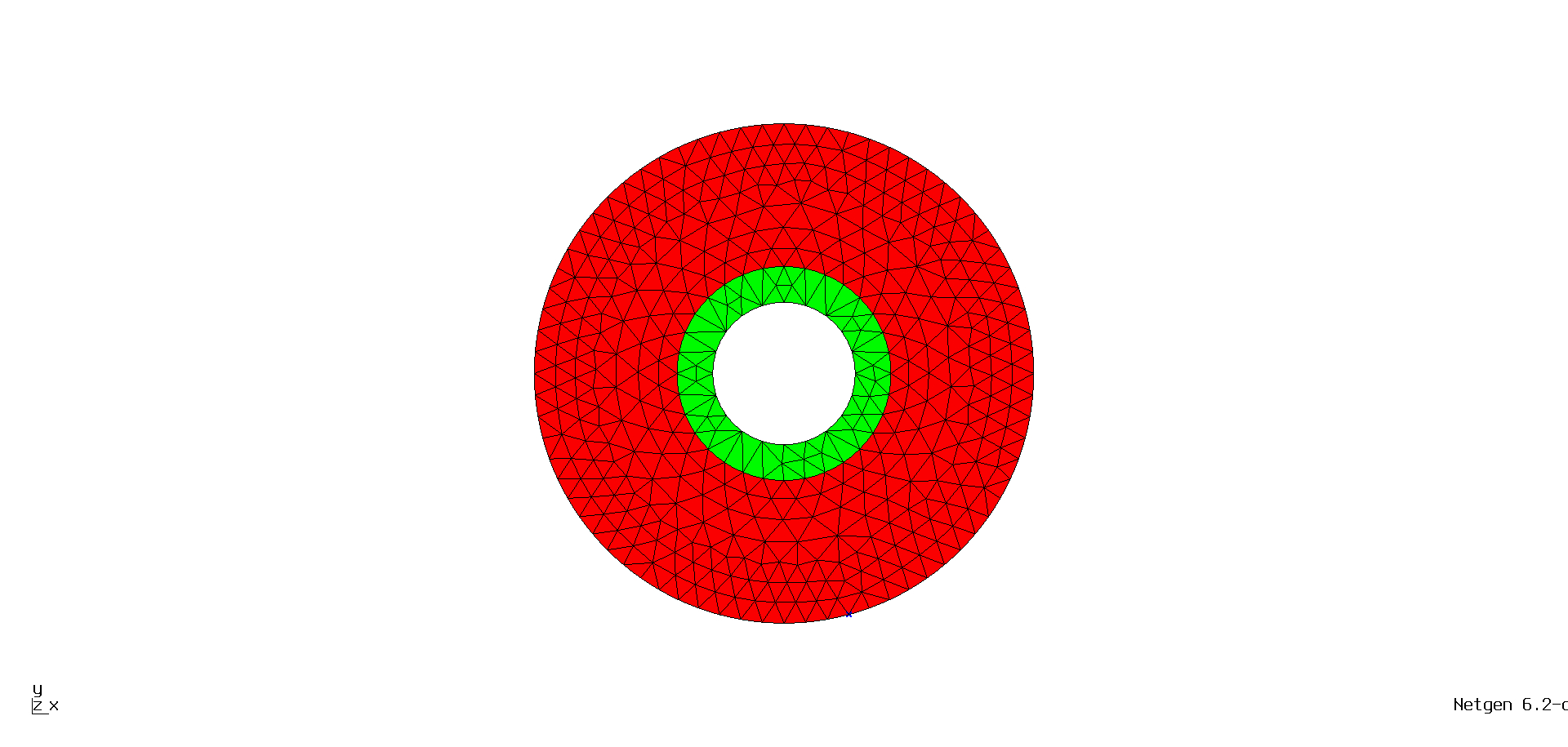}
\caption{Isotropic case.}
\label{fig:mesh_aniso}
\end{subfigure}
\caption{The computational domain with a coarse mesh for two geometries. The damping layer is color highlighted in red.}
\label{fig:meshes}
\end{figure}

\begin{figure}
\centering
\includegraphics[clip, trim={3cm 0cm 0cm 0cm}, scale=0.28]{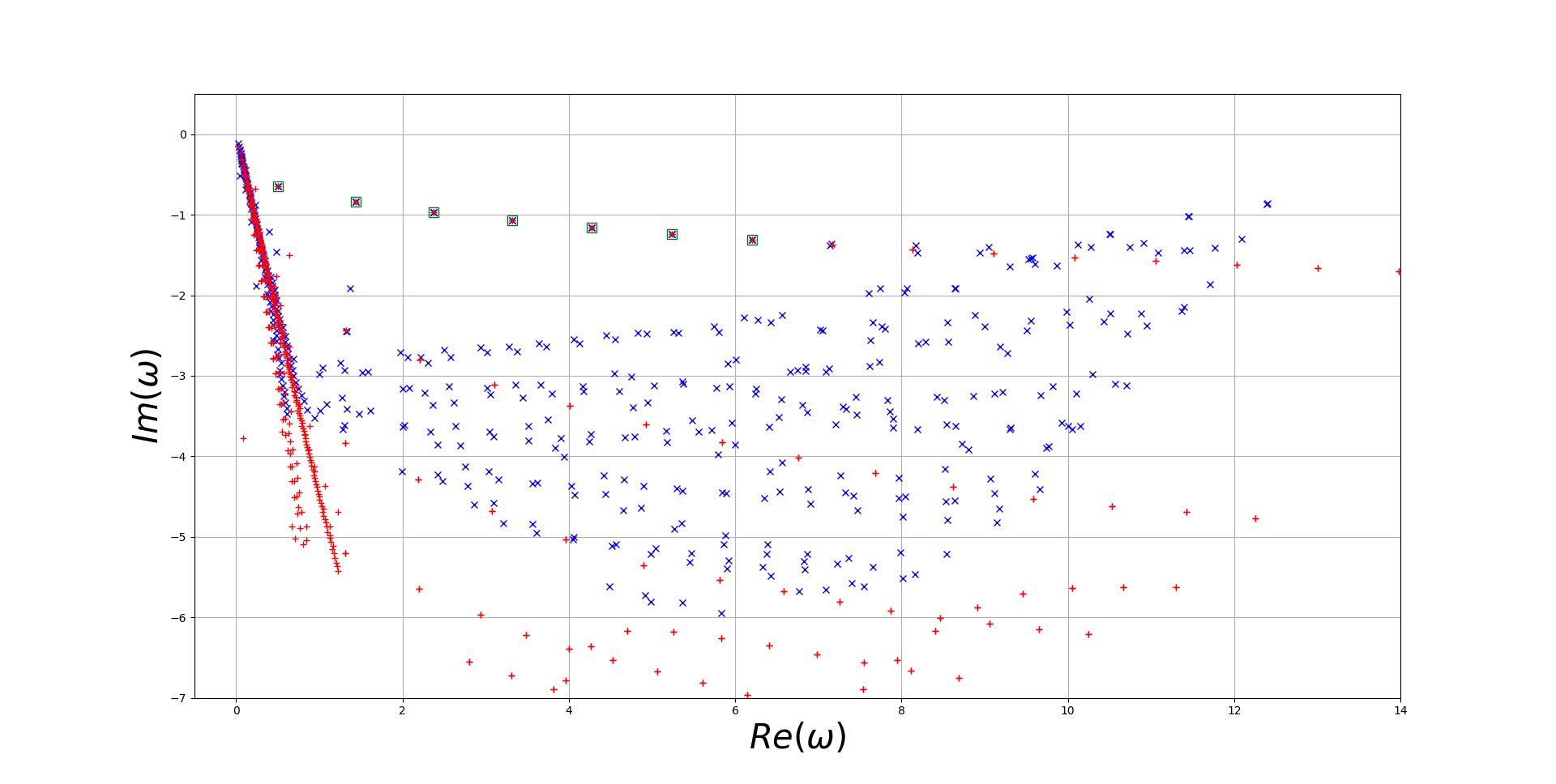}
\caption{Computed spectrum for the anisotropic geometry (blue x), computed spectrum for the isotropic geometry (red +), first seven semi-analytically computed resonances (green squares).
The parameters are $r_0=1.5$, layer width $L=2$, maximal mesh size $h=0.1$, uniform polynomial degree $p=6$, an affine scaling profile $\tilde\alpha(r)=1-r_0/r$ for $r>r_0$ and damping strength $\gamma=8i$.}
\label{fig:spectrum}
\end{figure}

In this section we verify our theoretical results with computational experiments.
Note that our preceding theory is formulated only for the 3D case with Dirichlet boundary conditions.
However, the analysis of variations such as the 2D Neumann case only requires minor adaptations.
To construct an example with available reference values we start with an isotropic problem.
We consider the 2D Neumann problem on $B_1^c$ with $\varsigma=I_{2\times 2}$.
For this setting it is well known that the resonances are the roots of the first derivatives of the Hankel functions of the first kind $(H^{(1)}_n)'$.
We compute these roots with the python routine cxroots as a set of (semi-)analytical reference values.
Now let $\Omega^c$ be the ellipse with principle axis $(-0.5,0.5)\times\{0\}$ and $\{0\}\times(-1,1)$.
Then by means of a domain transformation it is easy to see that the anisotropic 2D Neumann problem on $\Omega$ with $\varsigma=\bpm 0.25&0\\0&1\epm$ has the same resonances as the isotropic 2D Neumann problem on $B_1^c$.
Thus we apply the complex scaling technique to this anisotropic problem, truncate the domain and discretize with the finite element method.
Subsequently we compute the eigenvalues of this linear matrix eigenvalue problem with the Arnoldi algorithm.
To this end we use the finite element library NGSolve.
To have a second set of numerical reference values we also compute the resonances of the original isotropic problem in the same way.
In Figure~\ref{fig:meshes} we see the geometry and coarse meshes for the two case, and in Figure~\ref{fig:spectrum} we compare the computed spectrum with the reference values.
We observe in Figure~\ref{fig:spectrum} a good consistency between both computations and the semi-analytical reference values for the first few resonances, which confirms our theoretical results.
However, for the anisotropic case the pollution by spurious resonances is more invasive than for the isotropic case.
Note that these so-called spurious resonances are unconverged eigenvalues caused by a poor approximation of the resolvent at the given part of the spectrum, and they move to infinity for increased layer width/decreased mesh size (for more details on studies of spurious resonance see~\cite{NannenWess:18}).

\section{Conclusion}\label{sec:conclusion}

In this article we argued that \emph{radial} complex scaling/perfectly matched layer methods can be succesfully applied to anisotropic wave equations and do not suffer from the drawbacks of cartesian scalings~\cite{BecacheFauqueuxJoly:03}.
We presented heuristic reasons that this holds for any kind of wave equation (scalar, electromagnetic, elastic, \dots).
We analyzed the time-harmonic scalar case and its approximation in detail and showed that the techniques for the isotropic case required only little adaptations to treat the anisotropic case.
We presented computational examples to verify our theoretical results.
In comparison with the isotropic case it arose a new condition, which involves the CS/PML parameters and the degree of anisotropicy.
For time-dependent problems this implies that the convenient frequency dependency of the scaling may lead to instabilities.
We proposed an alternative frequency dependency as a remedy.

\bibliographystyle{abbrv}
\bibliography{short_biblio}

\begin{thebibliography}{10}

\bibitem{BecacheBonnetBDLegendre:06}
E.~B\'ecache, A.-S. {Bonnet-Ben Dhia}, and G.~Legendre.
\newblock Perfectly matched layers for time-harmonic acoustics in the presence
  of a uniform flow.
\newblock {\em SIAM J. Numer. Anal.}, 44(3):1191--1217, 2006.

\bibitem{BecacheFauqueuxJoly:03}
E.~B\'ecache, S.~Fauqueux, and P.~Joly.
\newblock Stability of perfectly matched layers, group velocities and
  anisotropic waves.
\newblock {\em J. Comput. Phys.}, 188(2):399--433, 2003.

\bibitem{BecFliKacKaz:21}
E.~B\'ecache, S.~Fliss, M.~Kachanovska, and M.~Kazakova.
\newblock On a surprising instability result of perfectly matched layers for
  {M}axwell's equations in 3d media with diagonal anisotropy.
\newblock {\em Comptes Rendus - Mathématique}, to appear.

\bibitem{BecacheKachanovska:17}
E.~B\'ecache and M.~Kachanovska.
\newblock Stable perfectly matched layers for a class of anisotropic dispersive
  models. {P}art {I}: necessary and sufficient conditions of stability.
\newblock {\em ESAIM: M2AN}, 51(6):2399--2434, 2017.

\bibitem{BecacheKachanovska:20}
E.~B\'ecache and M.~Kachanovska.
\newblock {Stability and Convergence Analysis of Time-domain Perfectly Matched
  Layers for The Wave Equation in Waveguides}.
\newblock Preprint, 2020.
\newblock hal-02536375.

\bibitem{Berenger:94}
J.-P. B{\'e}renger.
\newblock A perfectly matched layer for the absorption of electromagnetic
  waves.
\newblock {\em J. Comput. Phys.}, 114(2):185--200, 1994.

\bibitem{BermudezHervellaNPrietoRodriguez:08}
A.~Berm{\'u}dez, L.~Hervella-Nieto, A.~Prieto, and R.~Rodr{\'{\i}}guez.
\newblock An exact bounded perfectly matched layer for time-harmonic scattering
  problems.
\newblock {\em SIAM J. Sci. Comput.}, 30(1):312--338, 2007/08.

\bibitem{BonnetBDChambeyronLegendre:14}
A.-S. {Bonnet-Ben Dhia}, C.~Chambeyron, and G.~Legendre.
\newblock On the use of perfectly matched layers in the presence of long or
  backward guided elastic waves.
\newblock {\em Wave Motion}, 51(2):266--283, 2014.

\bibitem{BonnetBDetal:20}
A.-S. {Bonnet-Ben Dhia}, S.~N. Chandler-Wilde, S.~Fliss, C.~Hazard, K.-M.
  Perfekt, and Y.~Tjandrawidjaja.
\newblock The complex-scaled half-space matching method.
\newblock Preprint, 2020.
\newblock arxiv:2012.10721.

\bibitem{BonnetBDFlissTonnoir:18}
A.-S. {Bonnet-Ben Dhia}, S.~Fliss, and A.~Tonnoir.
\newblock The halfspace matching method: A new method to solve scattering
  problems in infinite media.
\newblock {\em Journal of Computational and Applied Mathematics}, 338:44--68,
  2018.

\bibitem{BramblePasciakTrenev:10}
J.~H. Bramble, J.~E. Pasciak, and D.~Trenev.
\newblock Analysis of a finite {PML} approximation to the three dimensional
  elastic wave scattering problem.
\newblock {\em Math. Comp.}, 79(272):2079--2101, 2010.

\bibitem{CassierJolyKachanovska:17}
M.~Cassier, P.~Joly, and M.~Kachanovska.
\newblock Mathematical models for dispersive electromagnetic waves: An
  overview.
\newblock {\em Computers and Mathematics with Applications}, 74
  (11):2792--2830, 2017.

\bibitem{ChewWeedon:94}
W.~C. Chew and W.~H. Weedon.
\newblock A 3{D} perfectly matched medium from modified {M}axwell's equations
  with stretched coordinates.
\newblock {\em Microwave Optical Tech. Letters}, 7:590--604, 1994.

\bibitem{CollinoMonk:98a}
F.~Collino and P.~Monk.
\newblock The perfectly matched layer in curvilinear coordinates.
\newblock {\em SIAM J. Sci. Comput.}, 19(6):2061--2090 (electronic), 1998.

\bibitem{DiazJoly:06}
J.~Diaz and P.~Joly.
\newblock A time domain analysis of {PML} models in acoustics.
\newblock {\em Computer Methods in Applied Mechanics and Engineering}, 195,
  29-32:3820--3853, 2006.

\bibitem{ErnGuermond:04}
A.~Ern and J.-L. Guermond.
\newblock {\em Theory and practice of finite elements}, volume 159 of {\em
  Applied Mathematical Sciences}.
\newblock Springer-Verlag, New York, 2004.

\bibitem{Halla:16}
M.~Halla.
\newblock Convergence of {H}ardy space infinite elements for {H}elmholtz
  scattering and resonance problems.
\newblock {\em SIAM J. Numer. Anal.}, 54(3):1385--1400, 2016.

\bibitem{Halla:19Diss}
M.~Halla.
\newblock {\em {A}nalysis of radial complex scaling methods for scalar
  resonance problems in open systems}.
\newblock PhD thesis, {T}echnische {U}niversit{\"a}t {W}ien, 2019.
\newblock https://repositum.tuwien.ac.at/urn:nbn:at:at-ubtuw:1-131893.

\bibitem{Halla:19Tcomp}
M.~Halla.
\newblock Galerkin approximation of holomorphic eigenvalue problems: weak
  {T}-coercivity and {T}-compatibility.
\newblock {\em Numer. Math.}, 2021.
\newblock https://doi.org/10.1007/s00211-021-01205-8.

\bibitem{Halla:20PML}
M.~Halla.
\newblock Analysis of radial complex scaling methods: scalar resonance
  problems.
\newblock {\em SIAM J. Numer. Anal.}, 2021 accepted.

\bibitem{HallaHohageNannenSchoeberl:16}
M.~Halla, T.~Hohage, L.~Nannen, and J.~Sch{\"o}berl.
\newblock Hardy space infinite elements for time harmonic wave equations with
  phase and group velocities of different signs.
\newblock {\em Numer. Math.}, 133(1):103--139, 2016.

\bibitem{HallaNannen:15}
M.~Halla and L.~Nannen.
\newblock Hardy space infinite elements for time-harmonic two-dimensional
  elastic waveguide problems.
\newblock {\em Wave Motion}, 59:94 -- 110, 2015.

\bibitem{HallaNannen:18}
M.~Halla and L.~Nannen.
\newblock Two scale hardy space infinite elements for scalar waveguide
  problems.
\newblock {\em Advances in Computational Mathematics}, 44(3):611--643, 2018.

\bibitem{HohageNannen:09}
T.~Hohage and L.~Nannen.
\newblock Hardy space infinite elements for scattering and resonance problems.
\newblock {\em SIAM J. Numer. Anal.}, 47(2):972--996, 2009.

\bibitem{HohageNannen:15}
T.~Hohage and L.~Nannen.
\newblock Convergence of infinite element methods for scalar waveguide
  problems.
\newblock {\em BIT Numerical Mathematics}, 55(1):215--254, 2015.

\bibitem{Joly:12}
P.~Joly.
\newblock An elementary introduction to the construction and the analysis of
  perfectly matched layers for time domain wave propagation.
\newblock {\em S$\vec{\rm e}$MA J.}, 57:5--48, 2012.

\bibitem{Moiseyev:98}
N.~Moiseyev.
\newblock Quantum theory of resonances: Calculating energies, width and
  cross-sections by complex scaling.
\newblock {\em Physics reports}, 302:211--293, 1998.

\bibitem{NannenWess:18}
L.~{Nannen} and M.~{Wess}.
\newblock {Computing scattering resonances using perfectly matched layers with
  frequency dependent scaling functions}.
\newblock {\em {BIT}}, 58(2):373--395, 2018.

\bibitem{NannenWess:19}
L.~Nannen and M.~Wess.
\newblock Complex scaled infinite elements for exterior helmholtz problems.
\newblock Technical report, 2019.

\bibitem{Simon:78}
B.~Simon.
\newblock Resonances and complex scaling: A rigorous overview.
\newblock {\em International Journal of Quantum Chemistry}, 14(4):529--542,
  1978.

\bibitem{Vavrycuk:07}
V.~Vavryčuk.
\newblock Asymptotic green's function in homogeneous anisotropic viscoelastic
  media.
\newblock {\em Proceedings of the Royal Society A: Mathematical, Physical and
  Engineering Sciences}, 463(2086):2689--2707, 2007.

\end{thebibliography}
\end{document}